\documentclass[letterpaper, 10 pt, conference]{ieeeconf}  % Comment this line out if you need a4paper

\IEEEoverridecommandlockouts                              % This command is only needed if
\usepackage{epsfig} % for postscript graphics files
\usepackage{mathptmx} % assumes new font selection scheme installed
\usepackage{times} % assumes new font selection scheme installed
\usepackage{amsmath} % assumes amsmath package installed
\usepackage{amssymb}  % assumes amsmath package installed
\usepackage{epstopdf}
\usepackage{color}
\usepackage{cite}

\newtheorem{alg}{Algorithm}
\newtheorem{assumption}{Assumption}

\title{\LARGE \bf
IPA in the Loop: Control Design  for Throughput  Regulation in Computer Processors
}

\author{X. Chen, Y. Wardi, and S. Yalamanchili
\thanks{Research supported in part by  the NSF  under Grant Number CNS-1239225.}
\thanks{School of Electrical and Computer Engineering, Georgia Institute of Technology, Atlanta, GA 30332, USA. Email:
xchen318@gatech.edu, ywardi@ece.gatech.edu, sudha@ece.gatech.edu.}
}

\begin{document}

\maketitle
\thispagestyle{empty}
\pagestyle{empty}

\begin{abstract}
A new technique for  performance regulation in event-driven systems, recently proposed
by the authors,  consists of
an adaptive-gain integral control. The gain is adjusted in the control loop  by a real-time estimation of the  derivative
of the plant-function with respect to the control input. This estimation is carried out  by Infinitesimal Perturbation
Analysis (IPA). The main motivation comes from applications to throughput regulation
 in computer processors,  where to-date, testing and assessment  of
the proposed control technique has been assessed by simulation. The purpose of this paper is to report on
its implementation on a machine, namely an Intel Haswell microprocessor, and compare its performance to that obtained
from cycle-level, full system simulation environment.
The intrinsic contribution of the paper to the Workshop on Discrete Event System is in describing the process of taking an IPA-based
design and simulation to a concrete implementation, thereby providing a bridge between theory and applications.

\end{abstract}

\section{Introduction}
One of the objectives   of systems' control is {\it performance regulation},
namely the output tracking of a given setpoint reference
 despite modeling uncertainties, time-varying  system's characteristics, noise, and other unpredictable
factors having the effects of system-disturbances. A commonly-practiced  way to achieve tracking is by a feedback control law that includes an integrator.
An integral control alone may have destabilizing effects on the closed-loop system, and hence the controller often
includes proportional and derivative elements as well thereby comprising
the well-known PID control \cite{Franklin15}.

Recently there has been a growing interest in performance regulation of event-driven systems, including  Discrete Event Dynamic Systems (DEDS) and Hybrid Systems (HS), and a control technique has been proposed  which   leverages on the
special structure of discrete-event dynamics \cite{Wardi16}.
The controller consists of a standalone  integrator with an adaptive gain, adjusted in real time as part of the control law. The rule for changing  the gain
is designed for stabilizing the closed-loop system as well as for simplicity
of implementation and robustness to computational and  measurement errors. Therefore it obviates the need for proportional and derivative elements, and can be implemented in real-time environments by approximating complicated computations by simpler ones. In other words,
the balance between precision and required computing efforts can be tilted in favor of simple, possibly imprecise  computations. A key feature of the control law is that it is based on the derivative
of the plant function,  namely the relation between  the system's  control parameter and its output, which is
computed or estimated by
Infinitesimal Perturbation Analysis (IPA). This will be explained  in detail in the following paragraphs.

IPA is a well-known and well-tested technique for computing sample-performance derivatives (gradients)
in DEDS, HS, and other event-driven systems  with respect to controlled variables; see \cite{Ho91,Cassandras99} for extensive presentations and surveys.
Its salient feature is in simple rules for tracking the propagations associated with a gradient along the
sample path of a system, by low-cost algorithms. However, this simplicity may come at the expense of statistical
unbiasedness of the IPA derivatives. In situations where IPA is biased, alternative perturbation-analysis techniques have been proposed, but they may require far-larger computing efforts than the basic IPA (see \cite{Ho91,Cassandras99}).
For the performance regulation technique described in this paper, it has been shown that IPA need not be unbiased and, as mentioned
earlier, its most important requirement is low computational complexity \cite{Wardi16}.

The control system we consider is depicted in Figure 1. Assuming discrete time and one-dimensional variables, $r$ is the
setpoint reference, $n=0,1,\ldots$ denotes the time counter, the control variable $u_n$ is the  input to the plant at
time $n$,
$y_n$ is the corresponding output, and $e_n:=r-y_n$ is the error signal at time $n$.
The control law is defined by the following equation,
\begin{equation}
u_n=u_{n-1}+A_n e_{n-1},
\end{equation}
and we recognize this as an adder, the discrete-time analogue of an integrator,
if the gain $A_n$ is a constant that does not depend on $n$. The plant is an event-driven  dynamical system whose
output $y_n$ is related to its input $u_n$ in a manner defined in the next paragraph, and denoted by the functional term
\begin{equation}
y_{n}=L_{n}(u_{n}),
\end{equation}
where $L_{n}:R\rightarrow R$ is a random function. Its IPA derivative $L_{n}^{\prime}(u_{n}):=\frac{\partial y_n}{\partial u_n}$ is used to
define the controller's gain $A_{n+1}$ via the equation
\begin{equation}
A_{n+1}=\Big(L_{n}^{\prime}(u_{n})\Big)^{-1},
\end{equation}
and the error signal is defined as
\begin{equation}
e_{n}=r-y_n.
\end{equation}
A recursive application of Eqs. (1)-(4) defines the closed-loop system.

As for the plant, it can have the following form. Consider a continuous-time  or discrete-time dynamical system whose input
is $u(t)\in R$, and its state is $z(t)\in R^q$ for some $q\geq 1$; the notation $t$ designates continuous time or discrete time. Let $g:R^q\rightarrow R$ be a function that is absolutely
integrable over finite-length intervals. Partition the time-axis $\{t\geq 0\}$ into
contiguous left-closed, right-open intervals, $C_1,C_2,\ldots$, called
{\it control cycles}. Suppose that the input to the latter dynamical system has  constant values
during each interval $C_{n}$, and it    can be changed only at the boundary of these intervals.  In the setting of
the system of Figure 1,  $u_{n}$ is the value of the input
$u(t)$ during $C_{n}$, and $y_n$ can be either
$
y_{n}:=\int_{C_{n}}g(z(t))dt
$ or $y_{n}:=\int_{C_{n}}g(z(t))dt/|C_{n}|$ where $|C_{n}|$ is the duration of $C_{n}$. In the case of discrete time,
 a sum-term of the form
$\sum_{k}g(z_{k})$ replaces the integral. We do not specify how to determine the control cycles  $C_{n}$, and they can have an a-priori
constant length, or their termination can be the result of certain events.

 For example,  let the plant consist of an M/D/1 queue, $u_n$ is the value of the
service time during $C_n$, and $y_n$ is the time-average of the sojourn times of all jobs arriving during $C_n$. IPA can be used
to compute the derivative $\frac{\partial y_n}{\partial u_n}$ via a well-known
formula \cite{Ho91}. Generalizing this example, suppose  that
the plant-system is a stochastic event-driven system (DEDS or HS), $u\in R$ is a control variable, assumed
to have a constant value ($u_{n}$) during $C_{n}$,  $y_{n}\in R$, is a random function of
$u_{n}$ as indicated by Eq. (2), and its derivative $L_{n}^{\prime}(u_{n})$ is computed by IPA. Later we will be concerned
with measurement and computational errors and hence modify Eqs. (2) and (3) accordingly.

The development of the proposed regulation technique has been motivated primarily by applications to
computer cores, especially regarding regulation of power and instruction-throughput by
adjusting the core's clock rate (frequency) \cite{Almoosa12a,Almoosa12b}. Concerning throughput, there are no prescriptive, let alone
analytic models for the frequency-to-throughput relationship, and a complicated, intractable queueing model has had
to be used for simulation.  Nonetheless a simple IPA algorithm has been developed and used to approximate the sample derivative for determining the integrator's
gain via Eq. (3). The regulation
technique was extensively tested on programs from an industry-based suite of benchmarks, Splash-2 \cite{Wu95}, using  a detailed simulation platform
for performance assessment of computer architectures, Manifold \cite{Wang14}.
We reported the results in \cite{Almoosa12b,Chen15,Wardi16}, and
deemed them encouraging and meriting a further exploration of the regulation technique.

In the context of IPA research, this regulation technique represents two new perspectives.  First, the traditional application of IPA
throughout its development has been to optimization, whereas here it is applied in a new way, namely to performance regulation.
Second, much of the research of IPA
has focused on its unbiasedness, whereas here, in contrast, the concern is with fast computation which may come
at the expense of accuracy and unbiasedness.   The main novelty of the paper as compared to
References  \cite{Almoosa12b,Chen15,Wardi16}
 is in the fact that it concerns not simulation but an actual implementation. In this we were facing
new challenges associated with real-time measurements,  computations, and control. Consequently we were unable to
control each core separately as in \cite{Almoosa12b,Chen15,Wardi16}, and hence applied the regulation
method to a processor   containing multiple cores. Furthermore, due to issues with real-time computation, we were forced
to take drastically cruder approximations to the IPA derivatives than in \cite{Almoosa12b,Chen15,Wardi16}, and in fact it seems that we drove the
degree of imprecision to the limit.
How this worked on application programs will be seen in the sequel. In any event, the work described here is, to our knowledge,  the first  implementation (beyond simulation) of IPA in a real-time control environment.

The rest of the paper is structured as follows. Section 2 summarizes  relevant convergence results of the regulation technique
in an abstract setting. Section 3 describes the system under study and its model, presents simulation results
on Manifold followed by implementation on a state-of-the-art computer processor, and compares the two.
Section 4 derives some lessons from these results and proposes directions for future research.

\vspace{.2in}
\begin{figure}[h]
\centering
\includegraphics[width=0.40\textwidth]{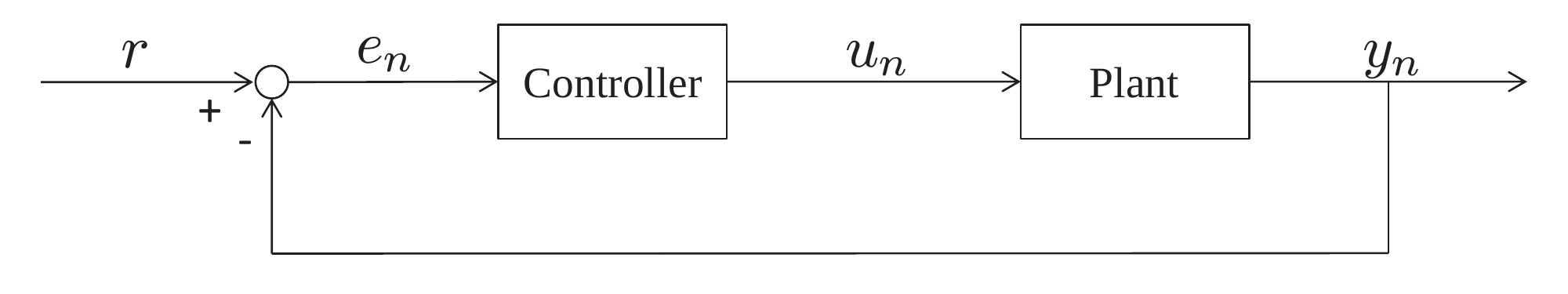}
% %\includegraphics[width=3.15in,angle=0]{Adaptive.jpg}
{\small \caption{Abstract regulation  system}}
\end{figure}

\section{Convergence  Results}
This section recounts established results concerning convergence of the regulation technique defined by recursive applications
of Eqs. (1) to (4), as summarized in Ref. \cite{Wardi16}. Ideally convergence means that
\begin{equation}
\lim_{n\rightarrow 0}e_{n}=0,
\end{equation}
hence $y_{n}\rightarrow r $ (see Figure 1). This can be achieved under suitable assumptions (mentioned below)  if the plant system is time invariant, and hence
 the function $L(u)=L_{n}(u)$  does not
depend on $n=1,\ldots$. In that case the control loop comprised of Equations (1)-(4) essentially implements Newton's method
for solving the equation
$r-L(u)=0$, for which there are well-known sufficient conditions for convergence. These include situations where the derivative term
$L^{\prime}(u_n)$ in Eq. (3) is computed in error, for which convergence in the sense of (5) is ascertained
under upper bounds on the magnitude of the error \cite{Lancaster66}.

If the system is time varying, Eq. (5) may not hold true, and in this case convergence can be characterized
by the equation
\begin{equation}
\limsup_{n\rightarrow\infty}|e_n|<\varepsilon,
\end{equation}
where $\varepsilon >0$ depends on a measure of the system's variability. To make matters concrete   let
$J:R\rightarrow R$ be a differentiable function,  and suppose that the term $y_{n}:=L_{n}(u_{n})$
in Eq. (2)  is a functional approximation of
$J(u_n)$. Assuming that $L_{n}(\cdot)$ is differentiable as well, we can view the term
$L_{n}^{\prime}(u_{n})$ as an approximation to $J^{\prime}(u_{n})$ in (3). However, for reasons that will be seen in the sequel,
we add another layer of approximation to $L_{n}^{\prime}(u_{n})$, denoted by $\zeta_{n}$, so that Eq. (3) computes the term
$L_{n}^{\prime}(u_{n})+\zeta_n$. Defining $\psi_{n}:=L_{n}(u_{n})-J(u_{n})$ and
$\phi_{n}=L_{n}^{\prime}(u_{n})+\zeta_{n}-J^{\prime}(u_{n})$, Eqs. (2) and (3) become
\begin{equation}
y_{n}=L_{n}(u_{n})=J(u_{n})+\psi_{n},
\end{equation}
and
\begin{equation}
A_{n+1}=\Big(L_{n}^{\prime}(u_{n})+\zeta_n\Big)^{-1}=\Big(J^{\prime}(u_{n})+\phi_{n}\Big)^{-1},
\end{equation}
respectively. The regulation technique now is defined by recursive applications of Eqs. (1),(7),(8),(4).

To analyze its convergence, suppose first, to simplify the discussion, that there exists a closed, finite-length interval,
 $I$, such that
every point $u_{n}$ computed by the regulation algorithm is contained in $I$; this assumption will be removed later.
Moreover, $I$ satisfies the following assumption.
\begin{assumption}
(i) The function $J$ and the functions $L_{n}$ are differentiable throughout $I$.
(ii) The function $J$ is either convex or concave, and monotone increasing or decreasing throughout $I$.
(iii) There exists $u\in I$ such that $J(u)=0$.
\end{assumption}

Various ways to relax this assumption will be discussed shortly. The following result was proved in \cite{Wardi16}.

{\it Proposition 2.3 and Lemma 2.2 in \cite{Wardi16}}:
For every $\varepsilon>0$, $\eta>0$,   and $M>1$,   there exist $\alpha\in(0,1)$,  $\delta>0$, and $\theta\in(0,1)$
 such that,  for every interval $I$ satisfying Assumption 1 and the following two additional conditions:  (i)
$
 \eta<=\min\{|J^{\prime}(u)|~:~u\in I\}$,  and (ii)
 $\max\{|J^{\prime}(u)|~:~u\in I\}\leq M\min\{|J^{\prime}(u)|~:~u\in I\}$:

 1). If for some $j\in\{1,2,\ldots\}$, and for $n=j,j+1,j+2$, $u_{n}\in I$,
 $|\phi_{n}|\leq\alpha|J^{\prime}(u_n)|$, and $|\psi_n|\leq\delta$,
 then
 \begin{equation}
|e_j|<\theta|e_{j-2}|.
\end{equation}

2). If for all $n=1,2,\ldots$, $u_{n}\in I$,
 $|\phi_{n}|\leq\alpha|J^{\prime}(u_n)|$, and $|\psi_n|\leq\delta$,
 then
\begin{equation}
\limsup_{n\rightarrow\infty}|e_n|<\varepsilon.
\end{equation}
\hfill$\Box$

In the context of the system considered in this paper, what we have in mind is a situation where the plant
is an event-driven system controlled by a real-valued parameter $u$, $J(u)$ is an expected-value function defined over a finite
horizon (hence not in steady state and possibly dependent on an initial condition),
$y_{n}=L_{n}(u_n)+\psi_{n}$ is a sample-based approximation  (possibly biased!) of $J(u_n)$ over the control cycle $C_{n}$, and $L_{n}^{\prime}(u)+\phi_n$ is a sample approximation of $J^{\prime}(u_n)$.

A few remarks concerning Assumption 1 are due.

1). The differentiability assumption is unnecessary, convexity of  $J(u)$ and almost-sure differentiability of
$L_{n}(u)$ at a given point $u$ suffice. These conditions often arise in the context of IPA. Under these weaker assumptions
 the proofs in
\cite{Wardi16} can be carried out in the context of convex analysis rather than differentiable calculus.

2). The condition that $u_{n}\in I$, $n=1,\ldots$, can be enforced in the case where $I$ is a constraint set for the
sequence $\{u_n\}$. In that case Eq. (1) would be replaced by
\begin{equation}
u_{n}=P_{I}(u_{n-1}+A_ne_{n-1})
\end{equation}
where $P_{I}(v)$ is the projection of $v\in R$ onto I, i.e., the point in $I$ closest to $v$.
The proof of convergence is unchanged.

3). Often the magnitude of the error terms $\psi_n$ and $\phi_n$ can be controlled by taking longer
control cycles, but there is no way to ensure the inequalities $|\phi_{n}|\leq\alpha|J^{\prime}(u_n)|$ and $|\psi_n|\leq\delta$
{\it for every} $n=1,2,\ldots$, which is stipulated as a condition for Eq. (10). Practically, however, with
 long-enough control cycles we
can expect those inequalities to hold for finite strings of $n\in\{k_1,\ldots,k_2\}$, thus guaranteeing the validity of
Eq. (9). If these strings are long enough,
Equation (9) would result in  $e_n$ approaching 0 at a geometric rate, then periodically jumping away due to the sporadic occurrence
of larger
errors, but again returning towards  0 rapidly, etc. This behavior has been observed in all of the examples where we tested
the regulation technique for a variety of event-driven systems \cite{Almoosa12b,Wardi16,Chen15,Seatzu14}.

4). Another source of the jitters described in the previous paragraph is the time-varying nature of the system.
This is particularly pronounced in the system tested in this paper, since the workload of programs processed by a core can
vary widely in unpredictable ways. Nonetheless we shall see that the regulation algorithm gives good results.

5).  It is possible to ascertain the assumptions underscoring the analysis in \cite{Wardi16} for simple systems
(e.g., tandem queues and some marked graphs), but  may be impossible to ascertain them for more complicated
 systems.  For instance, it can be impossible to prove differentiability or even convexity
of an expected-value function from characterizations  of its sample realizations, or bounds on the errors associated with
IPA. Moreover, some of these assumptions, including convexity or concavity, were made in \cite{Wardi16} in order to enable
an analysis but may not
be necessary. The aforementioned  convergence results serve to explain the behavior of the regulation method that was
observed in all of our former experiments \cite{Almoosa12b,Wardi16,Chen15,Seatzu14} as well as in those described
later in  this paper.

\section{Modeling, Simulation, and Implementation}

The system-architecture considered in this paper is based on an Out-of-Order (OOO) core technology
whereby instructions may complete execution in an order different from the
program order, hence the ``out of order'' designation. This enables instructions' execution
to be limited primarily by data dependency  and not by the order in which they appear in the program.
Data dependency arises when an instruction requires variables  that first must be computed by other instructions.
A detailed description of OOO architectures can be found in \cite{Hennessey12}, while a  high-level  description  is contained in \cite{Wardi16}.
Here we provide  an abstract functional and logical description, and refer the reader to \cite{Wardi16} for a
 more-detailed exposition.

The functionality of an OOO core is depicted in Figure 2.\footnote{Typically a core is dedicated to the processing
   of a program or a thread, namely a subprogram, as determined by the programmer or the operating system. In the forthcoming
   discussion we will use the term {\it program} to designate a thread as well.}
   Instructions are fetched sequentially  from memory and placed
in the instruction queue, where they are processed by functional units, or servers    in the parlance of queueing theory.
The queue is assumed to have unlimited storage and there is a server associated with each buffer.
The processing of an instruction starts as soon as it arrives {\it and} all of its required variables become available.
The instruction departs from the queue as soon as its execution is complete {\it and} the previous instruction
(according to the program order) departs.  In the parlance of computer architectures, an
instruction is said to be {\it committed} when it departs from the queue.  The instruction-throughput of the core is
defined and
measured by the average number of instructions committed per second.

\begin{figure}	[!t]
	\centering
	\includegraphics[width=3.0in]{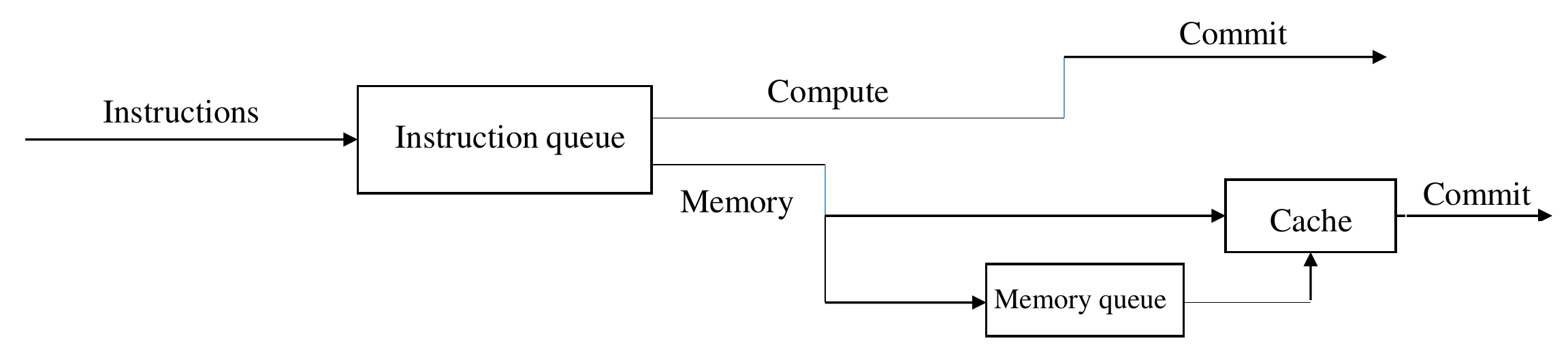}
	\caption{High-level instruction flow in OOO architectures}
	\label{OOO}
\end{figure}

Instructions generally are classified as {\it computational instructions} or {\it memory instructions}.
Access times of external, off-chip, memory
 instructions  are one-to-two orders of magnitude longer than those of
computational instructions.
Therefore most architectures make use of a hierarchical memory
arrangement where on-chip cache access takes less time than external memory such as DRAM. First the cache is searched,
and if the variable is found there then it is fetched and the instruction is completed. If the variable
is not stored in cache (a situation known as {\it cache miss}) then it is fetched from external memory
(typically DRAM) and placed in the cache, whence it is accessed and the instruction is  completed.
External  memory instructions can be thought of as being placed in
a finite-buffer, first-in-first-out queue,  designated as the {\it memory queue} in Figure 2. When this queue becomes full, the entire
memory access, including  cache, is stalled. Thus, there are three causes for an instruction to be stalled:
a computational or memory instruction waiting for variables computed by other instructions, a memory instruction waiting for the memory queue to
become non-full, and any instruction  waiting (after processing) for the previous instruction to depart from the queue.
We point out that instructions involving computation  and   $L_1$ cache-access are subjected
to the core's clock rate, while  memory instructions involving external memory, such as DRAM,  are not subjected to the same clock. This complicates
the application of IPA and may cause it to be biased.

A quantified discrete-event model of this process is presented in the appendix, and a more general description
can be found in \cite{Wardi16}, which also contains a detailed algorithm for the IPA derivative of the throughput as a function
of frequency.

\subsection{Manifold Simulation}
We use a  cycle-level, full system discrete event simulation platform for
multi-core architectures,  Manifold.
The simulated model consists of  a  16-core X86 processor die, where each core
 is in a separate clock domain and can control its own clock rate
 independently of other cores. For a detailed description of the Manifold simulation environment and its capabilities, please
 see Ref. \cite{Wang14}.

We simulated two application programs from the benchmark-suite Splash-2,
 Barnes and
Water-ns \cite{Wu95}. Barnes is a computation-intensive, memory-light application while Water-ns is  memory intensive.
 For each execution, all of the 16 cores of the processor run threads of
 the same benchmark concurrently while   each one of them is controlled
 separately.   The control cycle is
 set to $0.1$ ms for both Barnes and  Water-ns.
 The frequency range of the cores is set to $[0.5GHz,5.0GHz]$. We assume a continuous frequency range for the simulations, but later
 will consider a realistic, discrete range for the implementation described in the sequel.  The target instruction
 throughput   is set to the same value for each core
  for both Barnes and Water-ns, and we   experiment with the target throughput values of
   1,200 MIPS (Million Instructions per Second), 1,000 MIPS, and 800 MIPS. In terms of instructions per control
   cycle, these target values correspond to
   $0.12\times 10^6$, $0.1\times 10^6$, and $0.08\times 10^6$, respectively.
   The relationship between clock frequency and instruction throughput is determined by the Manifold processor model, but its IPA derivative
 was computed according to the high-level instruction flow described above and in the appendix.
 For
 each application run we present, in the following paragraphs,  the results for one of the 16 cores chosen at random.

\begin{figure}	[!t]
	\centering
	\includegraphics[width=3.5in]{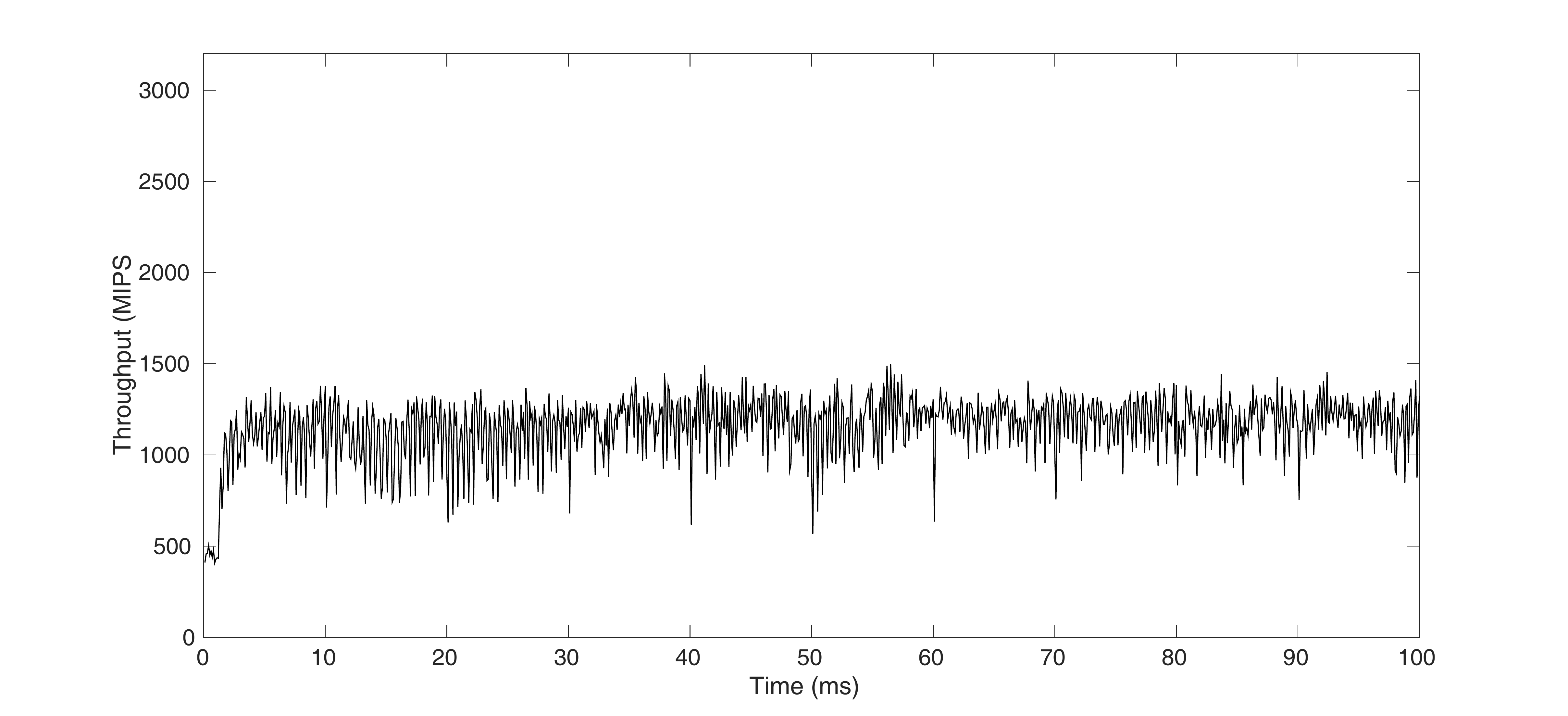}
	\caption{Manifold simulation, Barnes: throughput vs. time, target = 1,200 MIPS}
	\label{fig:MIPSManifoldBarnesCo10T1200}
\end{figure}

Consider first the target throughput of 1,200 MIPS.
The throughput simulation results for the Barnes benchmark  are shown in
Figure~\ref{fig:MIPSManifoldBarnesCo10T1200}, where the horizontal axis indicates time in ms and the vertical axis indicates
instruction throughput. We discern a rise of the throughput from its initial value (measured at      643.1 MIPS) towards
the target level of 1,200 MIPS, which it reaches for the first time in about
 1.5 ms, or 15 control cycles. Thereafter it oscillates about the target value,  which is not surprising  due to the unpredictable, rapidly-changing
program workload. The average throughput  computed over the time interval [$1.5ms$,$100ms$]
(soon after the throughput has reached the target value) is $1,157.4$ MIPS, which is 42.6 MIPS  off the target level of 1,200 MIPS.

\begin{figure}	[!t]
	\centering
	\includegraphics[width=3.5in]{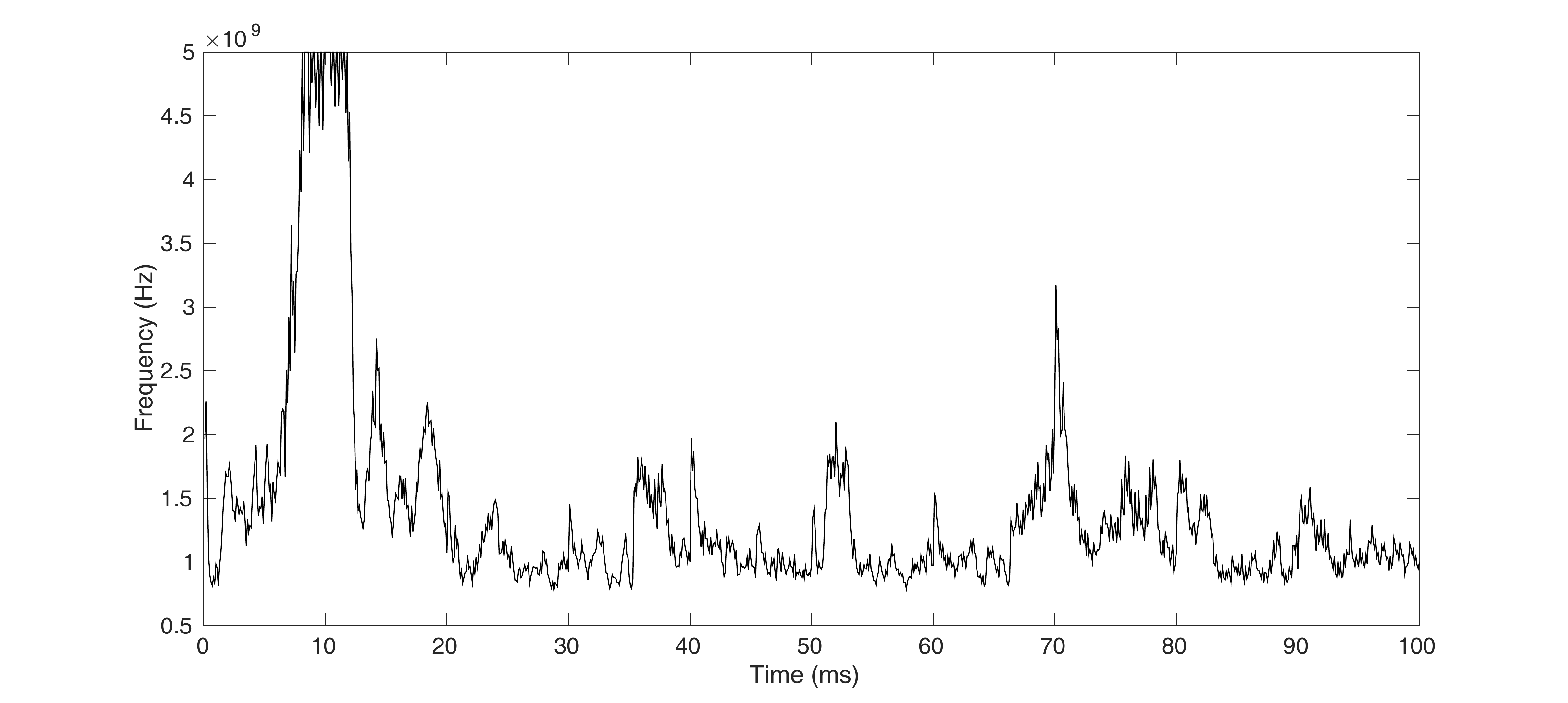}
	\caption{Manifold simulation, Barnes: clock frequency vs. time, target = 1,200 MIPS}
	\label{fig:FreqManifoldBarnesCo10Target1200}
\end{figure}

 Figure \ref{fig:FreqManifoldBarnesCo10Target1200} depicts  the graph of frequency vs. time (in ms), and it shows some saturation
 at its highest level of 5.0 Ghz, in the time-interval $[7ms,12ms]$. Saturation at the highest level can correspond to a negative offset of the average throughput from
 its target level, since it indicates that  the system may be unable  to raise  the throughput to a desired level.
 During the period of frequency saturation indicated in Figure 4, the throughput shown in Figure 3 it more jittery and
 sporadically attains slightly-lower values
 than after time 25ms. It also shows these characteristics between the end of the  saturation period and  time $25ms$. Therefore, the extent of
 the effects of the saturation on the
 aforementioned offset of 46.4 MIPS is not clear. Nonetheless we mention this point since it will be more pronounced in some of the results on which we report later.
  Also, we computed the average throughput in the intervals $[30ms,100ms]$ and $[50ms,100ms]$, after the jittery behavior
 of the throughput has somewhat subsided. The results are 1,192.6 MIPS and 1,192.9 MIPS, respectively, corresponding to offsets of
 7.4 MIPS and 7.1 MIPS from the target throughput of 1,200 MIPS.
 These results suggest that the frequency saturation plays some role in the larger, 46.4-MIPS
 offset  that was  computed over the interval
 $[1.5ms,100ms]$.

 For the target throughput of 1,000 MIPS,  the results (not shown due to space limitations) showed a rise in throughput from its initial value of
 420.5 MIPS to 1,000 MIPS in 2.1 ms, or 21 control cycles. The average throughput in the $[2.1ms,100ms]$ interval is 990.2 MIPS, corresponding to
 an offset of 9.8 MIPS of the throughput from its target value of 1,000 MIPS. The frequency saturated at its
 upper limit only at 5 isolated control cycles with minimal effects on the throughput.

 For the throughput target of 800 MIPS, the results show a rise in the throughput from its initial
 value
of 679.3 MIPS to 800 MIPS  in 1.9 ms, or 19 iterations.
The average throughput in the interval [1.9ms,100ms] was 839.6 MIPS, which
is 39.6 MIPS off the target value of 800 MIPS. There was a  considerable frequency
saturation at the lowest level of 0.5 GHz, which  explains the positive offset.

Returning to the results for the target level of 1,200 MIPS, we considered a way to reduce the  throughput oscillations and frequency saturation
by  scaling down the gain in Eq. (1). We did this by replacing
Eq.
(1) by the following equation,
\begin{equation}
u_n=u_{n-1}+\xi A_n e_{n-1},
\end{equation}
for a suitably-chosen constant $\xi\in(0,1)$. After some experimentation on various benchmarks (excluding those tested here) we chose
$\xi=0.2$. The resulting frequencies  did not saturate
throughout the program's run, and  yielded an average throughput of 1,198,5 MIPS, which is 1.5 MIPS off the target level of 1,200 MIPS.
Though working well for this example, this technique may be problematic when used with an implementation rather than simulation, as will be
discussed in the next subsection.

\begin{figure}	[!t]
	\centering
	\includegraphics[width=3.5in]{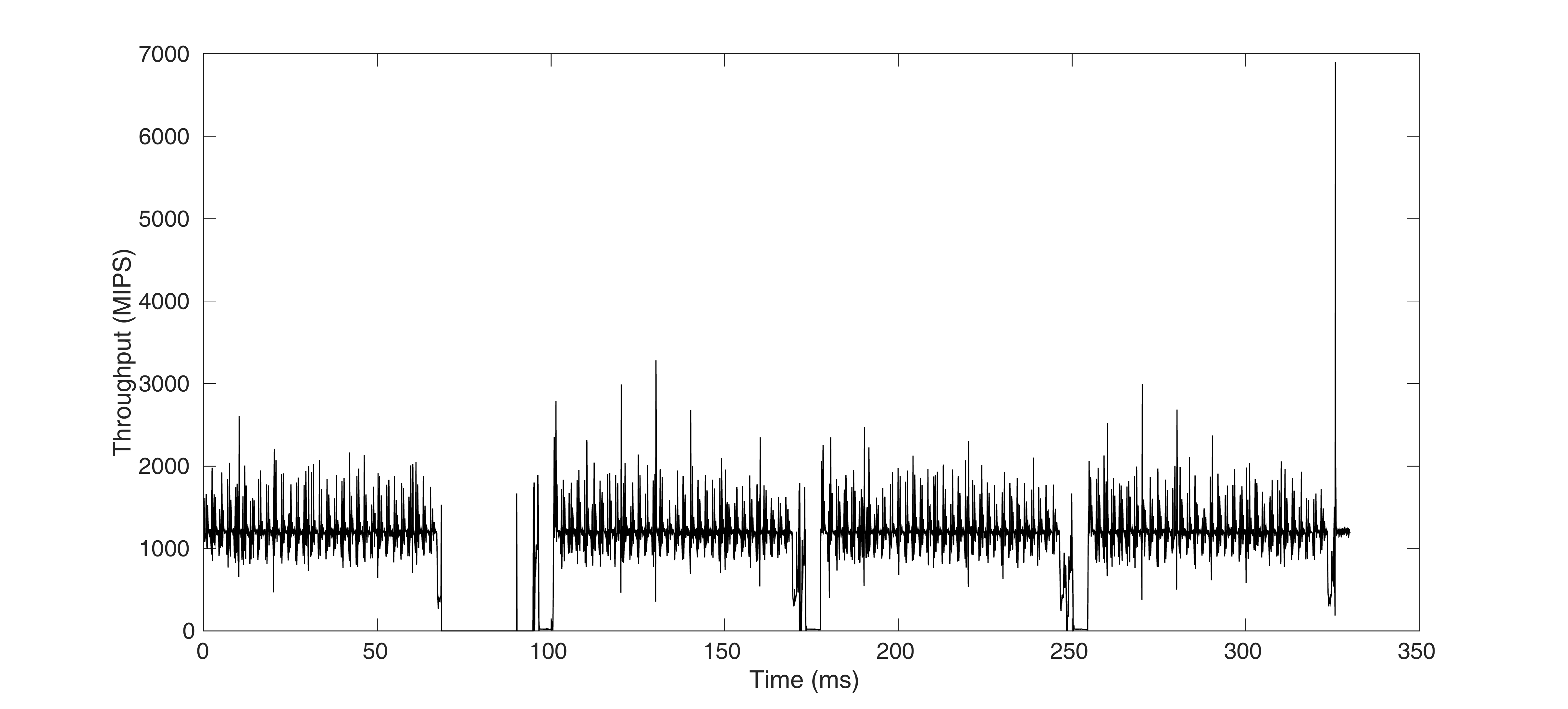}
	\caption{Manifold simulation, Water-ns: throughput vs. time, target = 1,200 MIPS}
	\label{fig:MIPSManifoldWaternsCo10Target1200}
\end{figure}

\begin{figure}	[!t]
	\centering
	\includegraphics[width=3.5in]{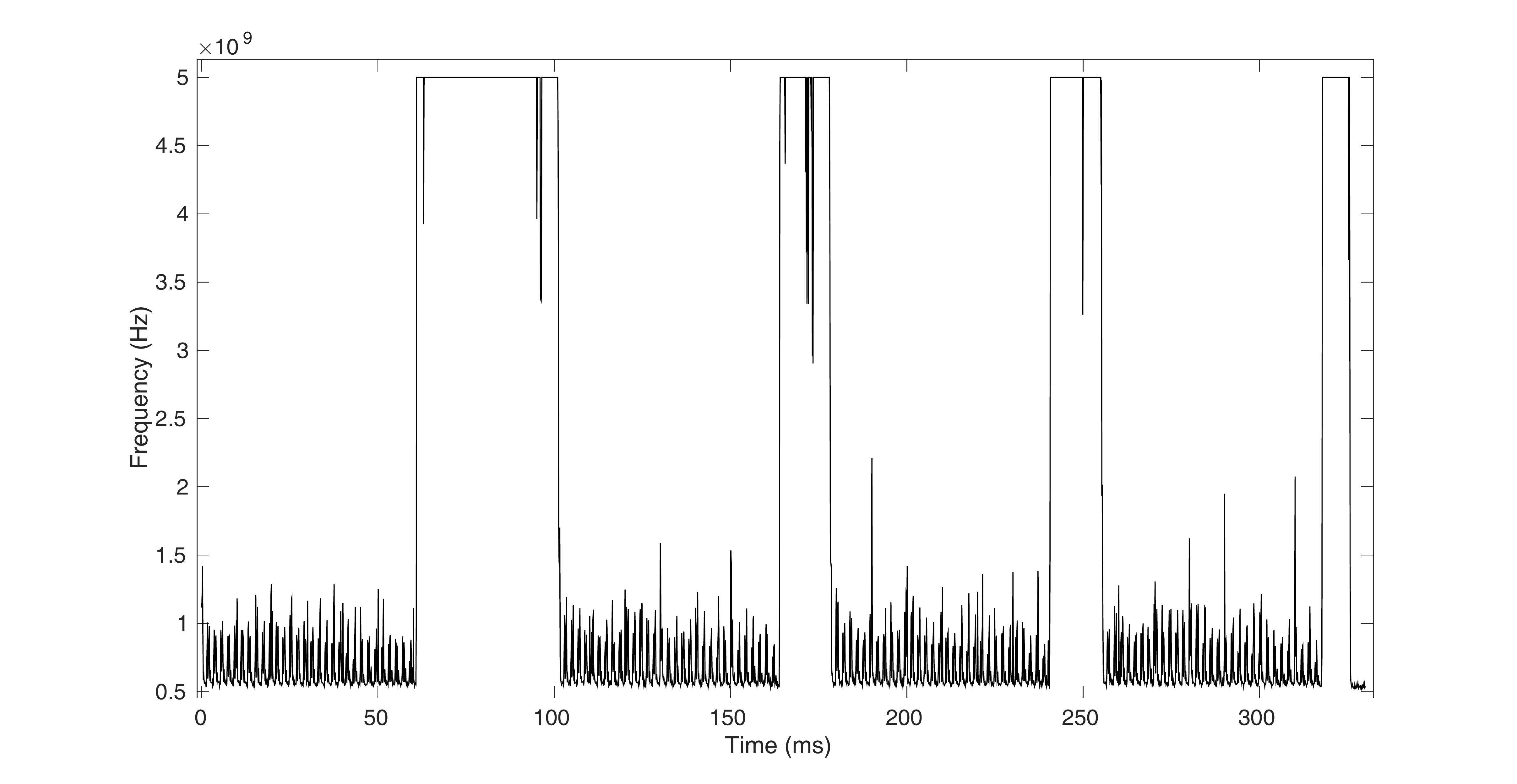}
	\caption{Manifold simulation, Water-ns: clock frequency vs. time, target = 1,200 MIPS}
	\label{fig:FreqManifoldWaternsCo10Target1200}
\end{figure}

For Water-ns, consider first the throughput target of 1,200 MIPS. Simulation results of throughput and frequency    are shown
in Figure~\ref{fig:MIPSManifoldWaternsCo10Target1200} and Figure \ref{fig:FreqManifoldWaternsCo10Target1200}, respectively.
We notice greater fluctuations and more saturation than for Barnes. In particular, Figure \ref{fig:FreqManifoldWaternsCo10Target1200}
shows three distinct periods of frequency saturations at its upper limit of 5.0 GHz, and Figure \ref{fig:MIPSManifoldWaternsCo10Target1200}
shows very low throughput during these periods. To explain this, recall that Water-ns is a memory-heavy
program, and execution times of memory
instructions are longer (typically by one or two orders of magnitude) than computational instructions. During those periods the instructions
of Water-ns mainly concern memory access, which are low-throughput instructions. The controller is applying its highest frequency in order
to push the throughput to its target value, but that frequency is not high
enough to have much effect. This is  why the periods of high-limit
frequency  saturation are
characterized by very low throughput. This has a pronounced affect of lowering the average frequency measured during the program's execution.  In fact, the  throughput obtained from the 
simulation rises from its initial value of 429 to its target level of 1,200 MIPS in about 1.8 ms, or
18 control cycles (this rise is not evident from Figure \ref{fig:MIPSManifoldWaternsCo10Target1200} due to
its insufficient granularity), and the 
 average throughout  from the time the target level is reached
($1.8ms$) to the end of the program-run ($333ms$)  is $1,126.8$ MIPS, which is 73.2 MIPS off the target level of 1,200 MIPS. 
Despite this offset, we observe 
that as soon as the program transitions from memory mode to computational mode, as indicated by the end of the saturation periods
in Figure  \ref{fig:FreqManifoldWaternsCo10Target1200}, the throughput returns quickly to about its target level, as can be
 seen in Figure  \ref{fig:MIPSManifoldWaternsCo10Target1200}.

 For the target throughput  of 1,000 MIPS, simulation results show the throughput
increasing from its initial value of
472.1 MIPS to the target level on 1,000 MIPS in 2.3 ms, or 23 control cycles.
There is considerable frequency saturation at the high limit of 5.0 GHz.
The average throughput in the interval $[2.3ms,330ms]$ is 947.8 MIPS, which means
an offset of 52.2 MIPS off the target throughput.

For the target throughput of 800 MIPS, simulation results indicated a rise of the throughput from its initial value of
443.3 to its target level in about 2.3 ms, or 23 control cycles. There is considerable saturation of the frequency at its lower level of
0.5 Ghz, and hence a positive offset between the computed average throughput and its target level.
Indeed, the average throughput in the interval $[2.3ms,330ms]$ is 862.6, MIPS, hence
 meaning an offset of 62.6 MIPS of the throughput from its
target value of 800 MIPS.

All of these results are summarized in Table I, showing the
offset (in MIPS) of average throughput from target throughput,  obtained from Manifold simulations of Barnes and Water-ns
 with throughput targets of 1,200, 1,00, and 800 MIPS.

Returning to the throughput target of 1,200 MIPS, an application of the modified algorithm with $\xi=0.2$ in Eq. (12)
yielded the average throughput of  1,143.6 MIPS which is 56.4 MIPS off the target level of 1,200 MIPS. This is a smaller offset than
the 73.2 MIPS obtained from the unmodified algorithm, and it is explained by the fact that there is still considerable,
though less frequency saturation than with the  unmodified algorithm.

%\begin{figure}	[!t]
%	\centering
%	\includegraphics[width=3.5in]{MIPSManifoldBarnesCo10Target800.pdf}
%	\caption{Manifold simulation, Barnes: throughput vs. time, target = 800 MIPS}
%	\label{MIPSManifoldBarnesCo10Target800}
%\end{figure}

\begin{table}[]
        \centering
        \caption{Manifold simulations: offset of average throughput from target levels}
        \label{my-label}
        \begin{tabular}{|l|l|l|l|}
                \hline
                Target Throughput (MIPS)& 1,200  & 1,000 & 800  \\ \hline  \hline
                Barnes            & -42.6 & -9.8 & 39.6 \\ \hline
                Water-ns          & -73.2 & -52.2   & 62.6 \\ \hline
        \end{tabular}
\end{table}

\subsection{Haswell implementation}
Haswell is Intel's fourth-generation core processor architecture fabricated in the 22nm process \cite{Hammarlund14}. Haswell
is comprised of multiple out-of-order execution cores designed for improved
power efficiency over prior generations.  The version used in our study has four cores residing on
each Haswell processor. Each core supports two  threads.  All the cores  execute at
the same frequency, the {\it processor frequency}.

Since the four cores operate at the same frequency, we cannot control each core separately
by the frequency. Instead, we consider the average  throughput among the active threads, which we call the {\it normalized processor throughput}, or {\it normalized throughput} in brief; it is the equivalent measure of the core's frequency in the Manifold
model described above.
We point out that typically the programmer and operating system
distribute the load among the various cores in a balanced way, and   in the system considered here there are 4
cores executing 8 threads, two threads to a core.

We implemented the controller by loading a C++ program to the processor via the PAPI
interface \cite{Browne00}.\footnote{Modern microprocessors include many hardware counters that record the occurrences of various events during program executions. Examples of such events include i) completion of the execution of an integer instruction, ii) a cache miss, or iii) an instruction that accesses memory.   The Performance Application Programming Interface (PAPI) is a publicly available software infrastructure for accessing these performance counters during program execution.}
Recall that the control algorithm
is based on  Eqs. (1),(2),(8),(4). However, the Haswell processor admits only a finite set of 16 frequencies,
 and  we have to modify Eq. (1) accordingly. This set of frequencies, denoted by
 $\Omega$, is (in GHz) $\Omega=\{0.8, 1.0, 1.1, 1.3, 1.5, 1.7, 1.8, 2.0, 2.2, 2.4, 2.5, 2.7, 2.9, 3.1,\\
  3.2, 3.4\}$.
 Denoting by $[u]$ the nearest point to $u\in R$ in the set $\Omega$ (the left point in case of a tie), we modify
  Eq. (1) by Eq. (13), below. The control algorithm is formalized as follows.

 {\it Notation}: $C_n$ - the $nth$ control cycle; $r$ - the target normalized throughput; $u_n$ -
 the processor frequency during $C_n$; $y_n$ - the resulting measured normalized throughput; $e_n:=r-y_n$.
\begin{alg}
The following steps are taken during $C_n$:
\begin{enumerate}
\item
At the start of $C_n$, set
\begin{equation}
u_n=\big[u_{n-1}+A_ne_{n-1}\big].
\end{equation}.
\item
During  $C_n$, measure $y_n$, and compute an approximation to the IPA derivative,
$\frac{\partial y_n}{\partial u_n}$, denoted by $\tilde{\lambda}_n$; note that
$\tilde{\lambda}_n=\frac{\partial y_n}{\partial u_n}+\zeta_n=L_{n}^{\prime}(u_n)+\zeta_n$ in Eq. (8).
\item
At the end of $C_n$, set $A_{n+1}=\big(\tilde{\lambda}_n\big)^{-1}$.
\item
At the end of $C_n$, compute $e_n:=r-y_n$.
\end{enumerate}
\end{alg}

\vspace{.15in}
The IPA algorithm used for the Manifold simulation is too complicated for a real-time implementation, and therefore
we explored approximations thereto with the objective of having them be as simple as possible. The simplest we could find
was
\begin{equation}
\tilde{\lambda}_n=\frac{y_n}{u_n},
\end{equation}
and we argue for that on the following grounds. The  frequency-to-throughput performance function $y_{n}:=L_{n}(u_n)$
is monotone increasing and has a linear component. As a matter of fact, if there are no external-memory instructions
then the throughput is linear in the core's frequency, as can be seen from the analysis in the appendix and in
\cite{Wardi16}. In that case $\tilde{\lambda}_n=\frac{\partial y_n}{\partial u_n}$. With  a mix of memory and computing instructions, the function $L_n(u_n)$ can be thought of as
having a linear component (with a time-varying slope) and a nonlinear component. We expect the error term
$\tilde{\lambda}_n-\frac{\partial y_n}{\partial u_n}$ to be larger, and consequently the regulation technique to work less well,
for intensive memory programs like Water-ns than for memory-light programs like Barnes. This is evident from
the testing we performed, whose results are presented  in the following paragraphs.

\begin{figure}	[!t]
	\centering
	\includegraphics[width=3.5in]{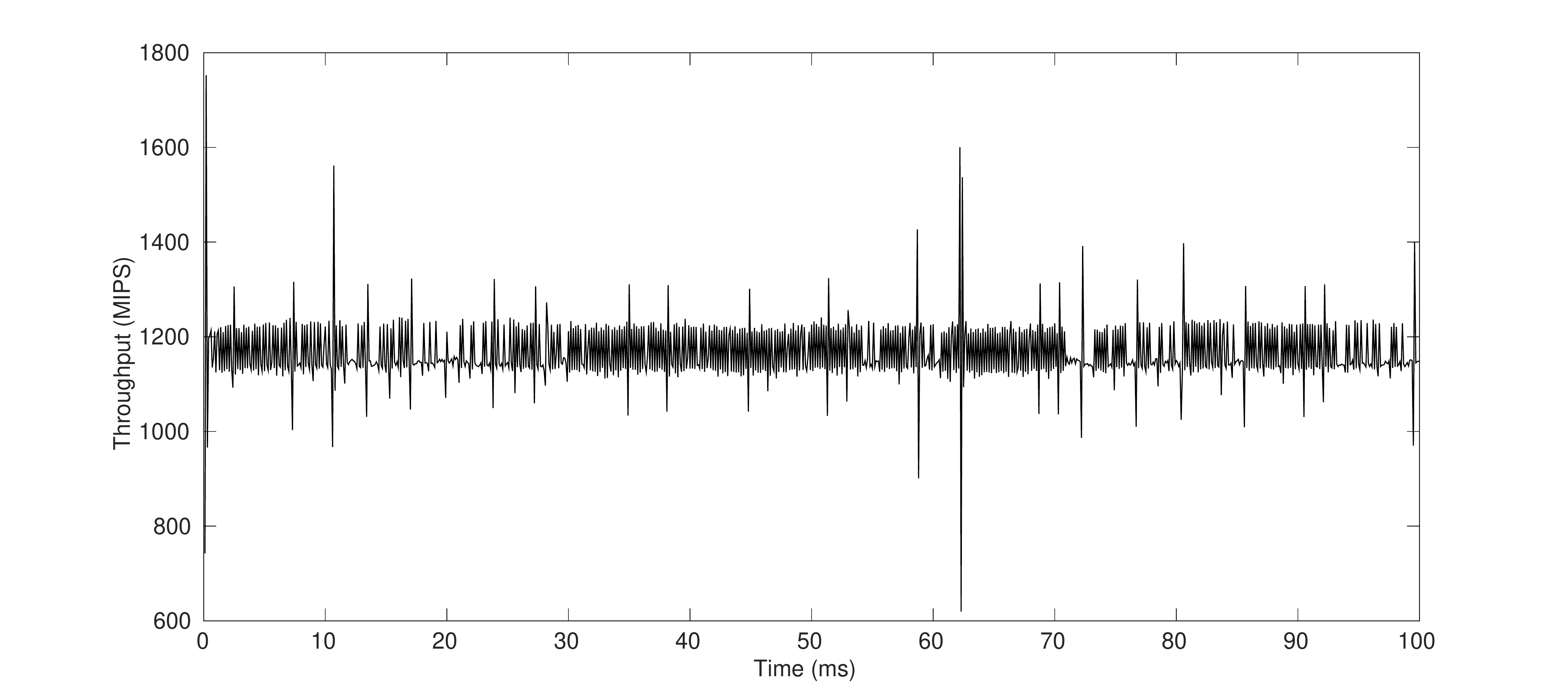}
	\caption{Haswell implementation, Barnes: throughput vs. time, target = 1,200 MIPS}
	\label{fig:MIPSHaswellBarnesCo10T1200}
\end{figure}

\begin{figure}	[!t]
	\centering
	\includegraphics[width=3.5in]{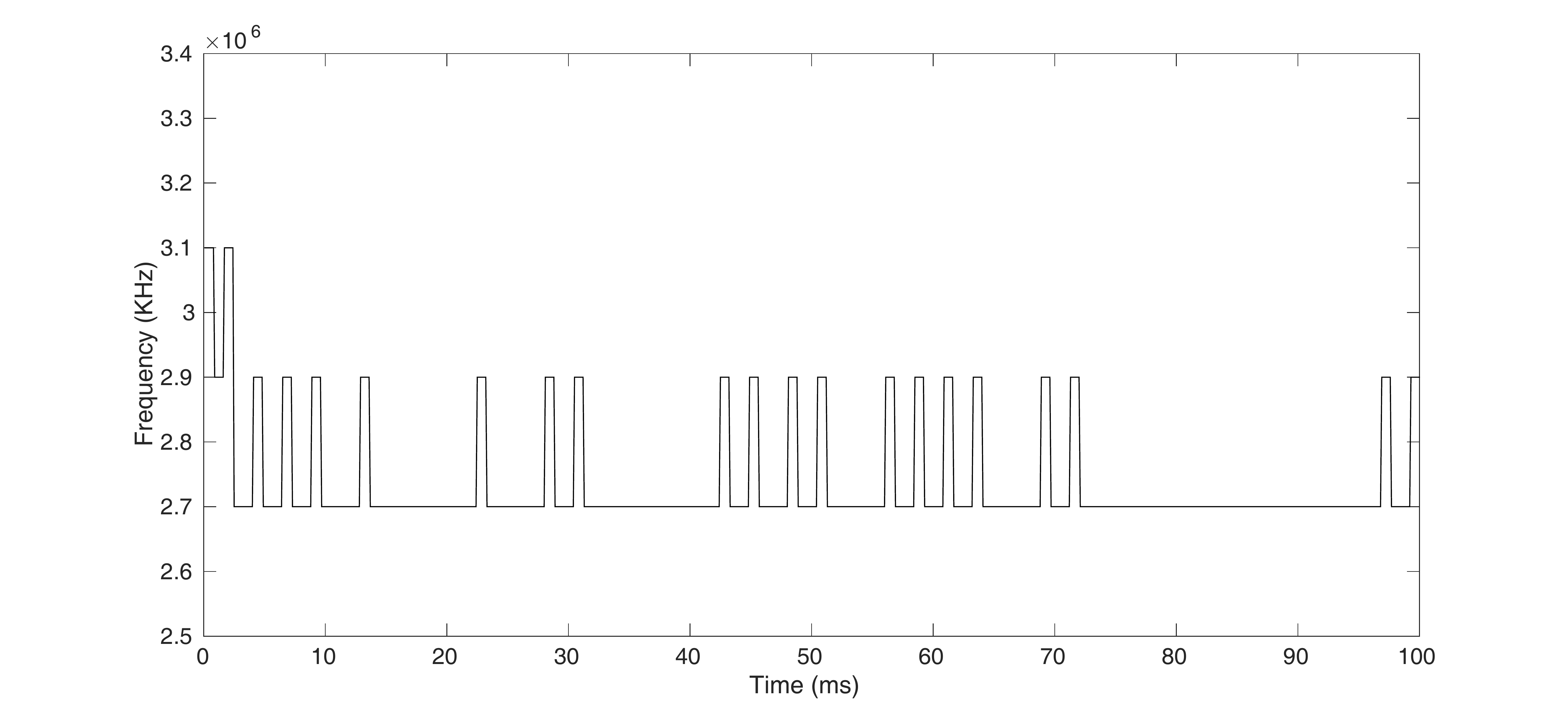}
	\caption{Haswell implementation, Barnes: clock frequency vs. time, target = 1,200 MIPS}
	\label{fig:FreqHaswelldBarnesCo10Target1200}
\end{figure}

Consider first the results obtained from the testing of Barnes.  For the throughput target of 1,200 MIPS, the results
are shown in Figure~\ref{fig:MIPSHaswellBarnesCo10T1200}, where the horizontal axis indicates time in ms and the vertical axis indicates instruction throughput.  The total run time is $100$ ms, and it corresponds to about 1,000 control cycles.
The throughput rises from an initial value of   $739.2$ MIPS to the target level of  1,200 MIPS in about $1.3$ ms, or  $13$ control
cycles. The average throughput computed over the time interval [$1.3ms$,$100ms$] (soon after the throughput has reached
the target value) is $1,166.5$ MIPS, which is 33.5  MIPS off the target level of 1,200 MIPS.
The  graph of the frequencies is shown in  Figure~\ref{fig:FreqHaswelldBarnesCo10Target1200}, and it indicates no saturation throughout the program.  We partly attribute the 33.5-gap to the quantization error due to the rounding off of the frequencies to their nearest values
in  $\Omega$, which is evident from Figure \ref{fig:FreqHaswelldBarnesCo10Target1200}.

For the target level of 1,000 MIPS, the throughput climbed from its initial value of 633.2 MIPS to its target level in
1.5 ms, or 15 control cycles. There was no frequency saturation, and the average throughput in the $[1.5ms,330ms]$ interval
is 990.6, which means an offset of 9.4 MIPS from the target level of 1,000 MIPS.

For the target level of 800 MIPS, the throughput climbs from its initial value
of
763.1 to the target level in 1.0 ms, or 10 control cycles. There was no frequency
saturation, and the average throughput is 829.7 MIPS, which is 29.7 MIPS off the target level of 800 MIPS.
Again, we attribute this gap to the saturation error.

Recall that we proposed a way to reduce the throughput oscillations and frequency saturation by modifying the control algorithm,
by replacing Eq. (1) by Eq. (12) with $\xi=0.2$. Although this  worked well for the Manifold simulation with the throughput target of
1,200 MIPS, it yielded poor results for the Haswell implementation.
 After a few iterations the
processor frequency was ``stuck'' at a value and did not move away from it. The reason is in the quantization error
inherent in the algorithm, which is due to the rounding off of the computed control variable in Eq. (13). The step size
for modifying the control variable  is insufficient to take that variable out of its current value.

\begin{figure}	[!t]
	\centering
	\includegraphics[width=3.5in]{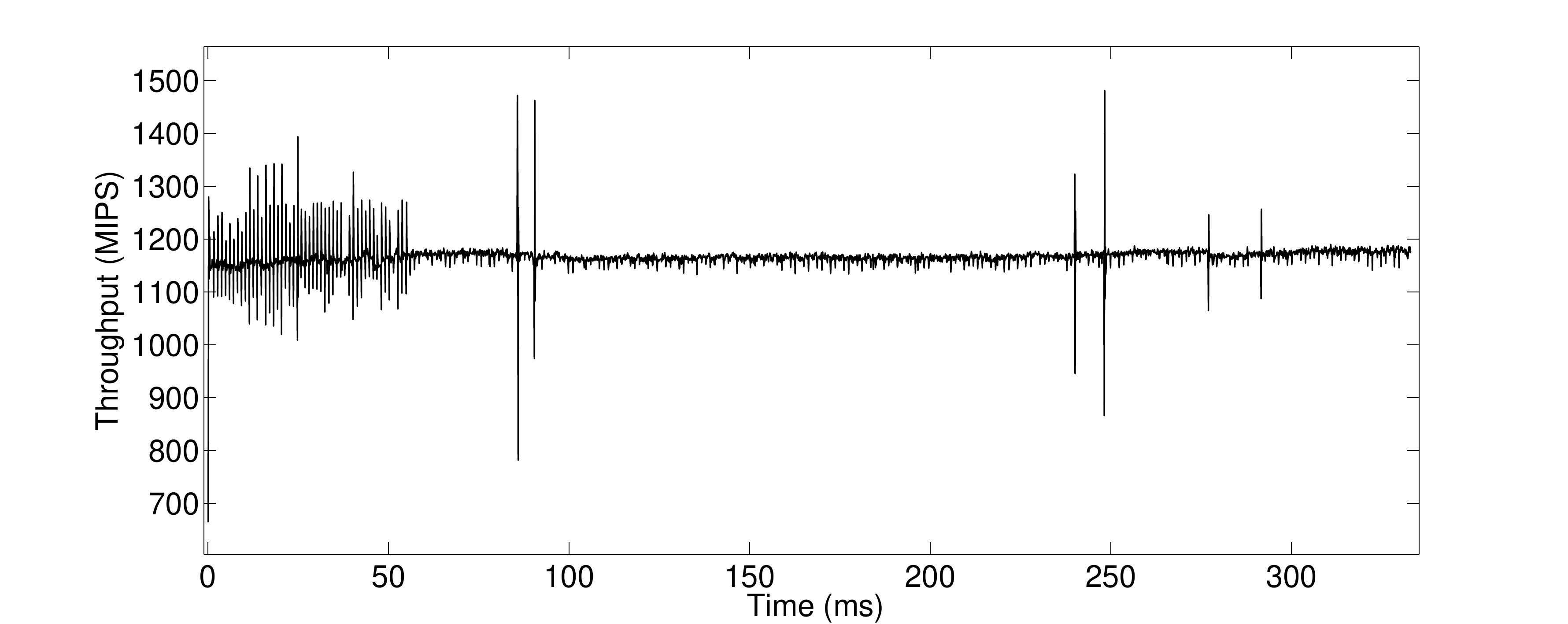}
	\caption{Haswell implementation, Water-ns: throughput vs. time, target = 1,200 MIPS}
	\label{fig:MIPSHaswellWaternsCo10Target1200}
\end{figure}

\begin{figure}	[!t]
	\centering
	\includegraphics[width=3.5in]{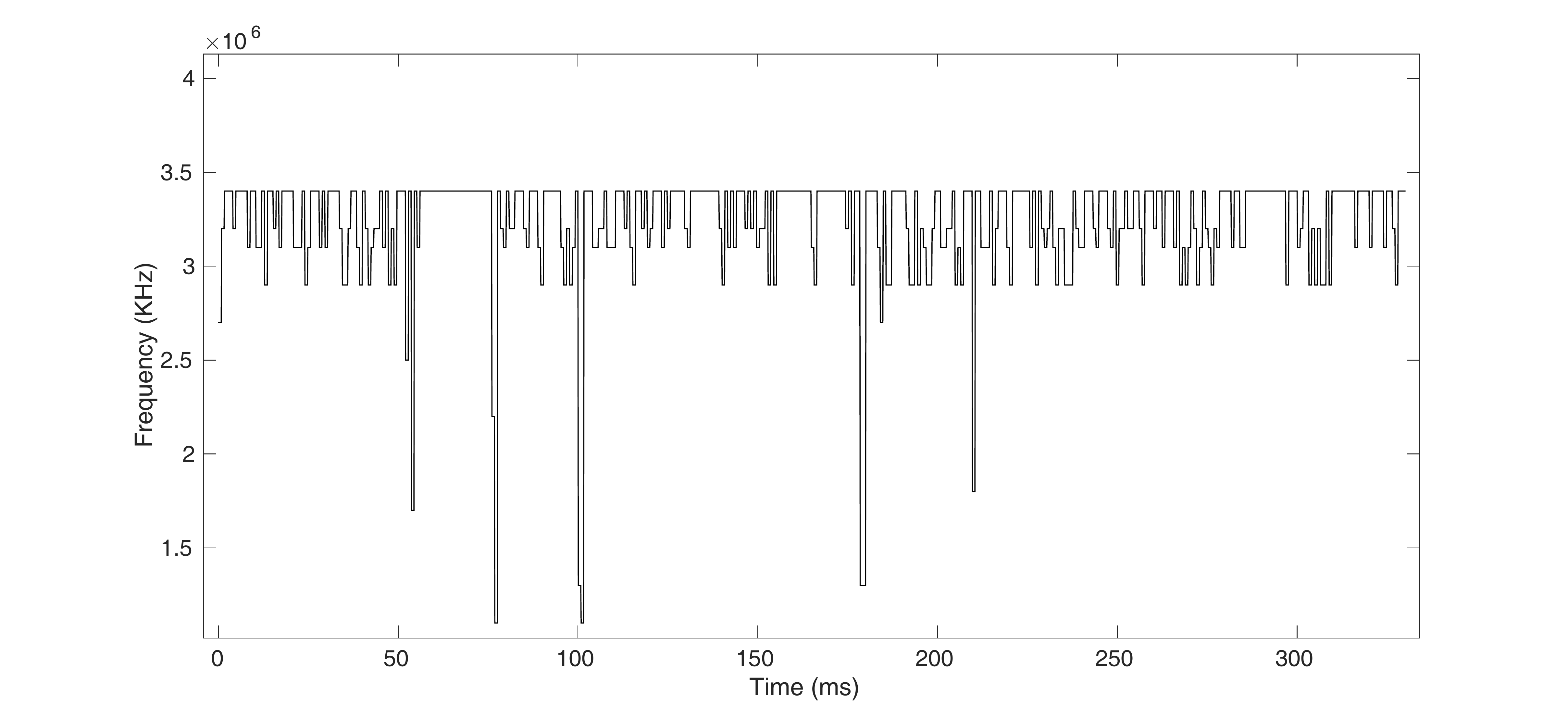}
	\caption{Haswell implementation, Barnes: clock frequency vs. time, target = 1,200 MIPS}
	\label{fig:FreqHaswellWaternsCo10Target1200}
\end{figure}

Regarding Water-ns, results for the throughput target of 1,200 MIPS are shown in  Figure~\ref{fig:MIPSHaswellWaternsCo10Target1200}.
The  throughput rises from
 an initial value of
 $683.1$ MIPS to its target value of 1,200 MIPS  in about $2.1$ ms, or  $21$ control cycles. The throughput oscillates at its upper level
 more than
 that obtained for Barnes, and this is due to the fact that Water-ns is a memory-heavy application. For the same reason,
 there is  considerable frequency saturation throughout the program as indicated in Figure
 \ref{fig:FreqHaswellWaternsCo10Target1200}. The obtained average throughput in the interval  [$2.1ms$, $333ms$] is $1143.2$ MIPS,
 which is 56.8 MIPS  off the target level of 1,200 MIPS.

 For the throughput level of 1,000 MIPS, the throughput rises from its initial value of 780 MIPS to its target level
 in 2.5 ms, or 25 control cycles. While there is some frequency saturation at the upper limit of the frequency range,
 the average throughput is 1,014.4 MIPS, which is 14.4 MIPS off its traget level.

Finally, for the target throughput of 800 MIPS, the throughput rises from its initial value of 698.3 to its target
level in 2.4 ms, or 24 control cycles. There are considerable frequency oscillations at the upper limit of the frequency range,
and the average throughput is   $836.3$, which is  36.3 MIPS  off the target level.

These results are summarized in Table II, showing the
offset (in MIPS) of average throughput from target throughput,  obtained from Haswell
implementation  of Barnes and Water-ns
 with throughput targets of 1,200, 1,00, and 800 MIPS.

\begin{table}[]
        \centering
        \caption{}
        \label{my-label}
        \begin{tabular}{|l|l|l|l|}
                \hline
                Target Throughput (MIPS) & 1,200  & 1,000 & 800  \\ \hline \hline
                Barnes            & -33.5 & -9.4 & 29.7 \\ \hline
                Water-ns          & -56.8 & 14.4   & 36.3 \\ \hline
        \end{tabular}
\end{table}

Comparing the data summarized in Table I and Table II, we that the regulation technique performs slightly  better  on the Haswell implementation platform
than on the Manifold simulation environment. The reason for this may be due to the fact
that in the simulation experiment we regulate the throughput of each core separately, while in the implementation we control 
the average throughput of all the cores in the processor.

\section{Conclusions}

This paper describes the testing of an IPA-based throughput regulation technique in multicore processors. The testing was
performed on both a simulation environment and an implementation platform. Despite crude approximations that have had to be made in
the implementation setting, the proposed technique performed slightly better than in the simulation setting. Future research will extend the
regulation method from a centralized control of a single processor to a distributed control of
networked systems.

\section{Appendix}
This section provides a quantitative description of the instruction-flow in the OOO-cache high-level model described at the
beginning of Section III.

Denote by $I_i$, $i=1,2,\ldots$, the instructions arriving at the instruction
queue in increasing order. Let $u$ denote the clock rate, or frequency, and let $\tau:=u^{-1}$ be the clock cycle. Denote by
$a_i(\tau)$ the arrival time of $I_i$ relative to the arrival time of $I_1$, namely $a_1(0):=0$,
and let $\xi_i$ be the clock counter at which $I_i$ arrives.
 Then,
$
a_i(\tau)=\xi_i\tau.
$
Denote by $\alpha_i(\tau)$   the time at which execution of $I_i$ starts, and let  $\beta_i(\tau)$ denote the time
at which execution of $I_i$ ends.

We next describe a way to compute  $\alpha_i(\tau)$.
consider first the case were $I_i$ is a computational instruction. If all of its required variables are available at its
arrival time then $\alpha_i(\tau)=a_i(\tau)+\tau$. On the other hand, if $I_i$ has to wait for such variables, let
$k(i)$ denote the index (counter) of the instruction last to provide such a variable, then
$\alpha_i(\tau)=\beta_{k(i)}(\tau)+\tau$. Next, if $I_i$ is a memory instruction, then $\alpha_i(\tau)$ is the time it
starts a cache access. If the memory queue is not full at time $a_i(\tau)$, then $\alpha_i(\tau)=a_i(\tau)+\tau$. On the other hand, if the memory queue is full at time $a_i(\tau)$, let $\ell(i)$ denote the index of the instruction at the head of the queue,
then, $\alpha_i(\tau)=\beta_{\ell(i)}(\tau)+\tau$.

To compute $\beta_i(\tau)$, consider first the case where $I_i$ is a computational instruction.
Let $\mu_i$ denote the number of clock cycles it takes to
execute $I_i$. Then,
$\beta_i(\tau)=\alpha_i(\tau)+\mu_i(\tau)$. On the other hand, if $I_i$ is a
memory instruction, let $\nu_i$ denote the number of clock cycles it takes to perform a
cache attempt. If the cache attempt is successful and the variable is found in cache, then $\beta_i(\tau)=\alpha_i(\tau)+\nu_i(\tau)$.
If the variable is not in cache, the instruction is directed to the
memory queue. Its transfer there involves a small number of clock cycles, $m_i$, hence its
arrives at the queue at time  $\alpha_i(\tau)+\nu_i\tau+m_i\tau$. The memory queue is a FIFO queue whose service
 time represents an external-memory access, which is independent of the core's clock. Denote by $S_i$ the sojourn time of $I_i$ at the memory
 queue. Then  $\beta_i(\tau)=\alpha_i(\tau)+\nu_i\tau+m_i\tau+S_i+\tau$.

 Finally, the departure time of $I_i$ from the instruction queue, denoted by $d_i(\tau)$, is
 $d_i(\tau)=\max\big\{\beta_i(\tau),d_{i-1}(\tau)\big\}+\tau$. Given a control cycle consisting of $N$ instructions, the throughput
 is defined as $N/d_{N}(\tau)$. Since $u=\tau^{-1}$, we can view the throughput as a function of $u$ and denote it by
 $y(u)$.
Its IPA derivative, $\frac{\partial y}{\partial u}$, can be computed by following the above dynamics of the instructions' flow.
This, and a more detailed discussion of the model, can be found in \cite{Wardi16}.

\end{document}